\documentclass[12pt]{elsarticle}
\usepackage[left=1in,right=1in, top=1.2in,bottom=1.2in]{geometry}
\usepackage{times}
\usepackage{amssymb,amsthm,latexsym,amsmath,,epsfig,pgf}
\usepackage{graphicx}
\usepackage{url}
\usepackage[T1]{fontenc}

\newtheorem{theorem}{Theorem}[section]
\newtheorem*{theorem A}{Theorem A}
\newtheorem*{theorem B}{N\"olker's Theorem}

\theoremstyle{remark}
\newtheorem{remark}{Remark}[section]
\theoremstyle{remark}

\begin{document}

\begin{frontmatter}
\title{Intersection Graphs of Graphs and Hypergraphs: A Survey}



\author[label1]{Ranjan N. Naik}
\address[label1]{\small Department of Mathematics,\\
Lincoln University, PA, USA\\

\vspace*{2.5ex} 
 {\normalfont rnaik@lincoln.edu}
 }
 
\begin{abstract}
The survey is devoted to the developmental milestones on the characterizations of intersection graphs of graphs and hypergraphs.  The theory of intersection graphs of graphs and hypergraphs has been a classical topic in the theory of special graphs. To conclude, at the end, we have listed some open problems posed by various authors whose work has contributed to this survey and also the new trends coming out of intersection graphs of hyeprgraphs.

\end{abstract}


\end{frontmatter}

\section{Introduction}
We follow the terminology of Berge, C. [3] and [4]. This survey does not address intersection graphs of other types of graphs such as interval graphs etc. An introduction of intersection graphs of interval graphs etc. are available in Pal [40]. A graph or a hypergraph is a pair (V, E), where V is the vertex set and E is the edge set a family of nonempty subsets of V. Two edges of a hypergraph are $l$ -intersecting if they share at least $l$ common vertices. This concept was studied in [6] and [18] by Bermond, Heydemann and Sotteau. A hypergraph H is called a k-uniform hypergraph if its edges have k number of vertices. A hypergraph is linear if any two edges have at most one common vertex. A 2-uniform linear hypergraph is called a graph or a linear graph. 

The intersection graph of a hypergraph H, denoted by L(H), is the graph whose vertex set is the edge set of H, and two vertices are adjacent in L(H) if the corresponding edges are adjacent or intersecting edges in H. Among the earliest works on intersection graphs of graphs is seen in Whitney's work [51] in 1934, Krausz's work [25] in 1943 and Hoffman's work [20] in 1964. Intersection graphs of graphs are also called line graphs. Harary [14] noticed that the intersection graph of a graph has been introduced under different names by many authors. Berge [4] called intersection graphs as representative graphs.The intersection graph of a hypergraph in [4] is a $l$-intersection graph with $l=1$, the adjacent edges may have one or more common vertices. 

We call $l$ the multiplicity of a hypergraph if any two edges of H intersect in at most $l$ vertices, $1  \leq  l $ $  \leq  k.$  We denote by $ I_{l} (k)$, the set of all graphs which are intersection graphs of some k-uniform hypergraphs with multiplicity $l$.Thus $I_{l} (k)$ is the set of graphs which are intersection graphs of some k-uniform linear hypergraphs. We write $I(k)$ for the set $ U_{(l \geq 0)} I_{l} (k)$. 

This survey is comprised of sections depending on the values of k and covers almost all the known results on the characterizations of intersection graphs of hypergraphs. Readers interested in proofs may refer to the actual papers cited in reference section.

\section{The characterization of intersection graphs of graphs.}

The classes $I_{1} (k)$ and $I (k)$ have been studied for a long time. The following well-known theorem on edge isomorphism on line graph of a graph is from Whitney.

\begin {theorem}$[51]$ If G and H are connected graphs and L(G) $\cong$L(H), then G $\cong$H unless one is $K_{3}$ and the other is $K_{1,3}$. Further, if the orders of G and H are greater than 4, then for any isomorphism $ \begin{matrix} f\colon & \mathbb{L(G)} \to \mathbb{L(H)} \end{matrix}$, there exists a unique isomorphism between G and H inducing $f$.
\end{theorem}

The other well-known characterization of line graph of a graph given below in theorem 2.2 from Krausz gives a global characterization of line graphs. This result is instrumental in the development of characterizations of intersection graphs of hypergraphs.

\begin{theorem}
$[25]$ A graph G is a line graph of a graph if the edges of G can be partitioned into complete subgraphs (cliques) in such a way that no vertex lies in more than two of the cliques and any two cliques have at most one common vertex.
\end{theorem}

Berge and Rado [5] and Gardner [12] obtained results for hypergraphs analogous in some sense to theorem 2.1. A relationship for hypergraphs which is stronger than isomorphism has been studied in [27], [47], and [49] by unifying the theorems of Krausz and Whitney.

The following theorem 2.3 of Rooij and Wilf describes in a structural criterion for a graph to be a line graph of a graph. A triangle T$(K_{3})$ of a graph G is said to be odd if there is a vertex of G adjacent to an odd number of vertices of T. A triangle T is said to be even if it is not odd. 

\begin{theorem}
$[42]$ A graph G is a line graph of a graph if G has no $K_{1,3}$, and if two odd triangles have a common edge, then the subgraph induced by their vertices is $K_{4}$.
\end{theorem}

A family of graphs M is said to be hereditary if $G \in M$ implies all the vertex induced subgraphs of G are also in M. The families $I(k)s$ are clearly hereditary families of graphs. Now let M be a hereditary family of graphs. If G does not belong to M, then clearly G is not an induced subgraph of any graph in M. A graph G does not belong to M is said to be a minimal forbidden graph for M if all vertex induced subgraphs of G are in M. Let $F(M)$ denote the family of all minimal forbidden graphs for M. Clearly, then M can be characterized by saying that $G \in M $ if and only if any graph of $F(M)$ is not an induced subgraph of G. We call such characterization the forbidden induced subgraph ($FIS$) characterization of the family of graphs. The well-known result of Beineke [2] and the unpublished work of Robertson N., stated below shows that $F(I_{1} (2))$ consists of the nine graphs  which cannot occur in line graphs of graphs. Beineke called line graphs derived graphs in [2].

\begin{theorem}
$[2]$ A graph G is an intersection graph of a graph if and only if none of the nine graphs given below in figure 1 from $[59]$ is an induced subgraph of G.
\end{theorem}
\begin{figure}[!htb]
    \center{\includegraphics[width=0.5\textwidth] {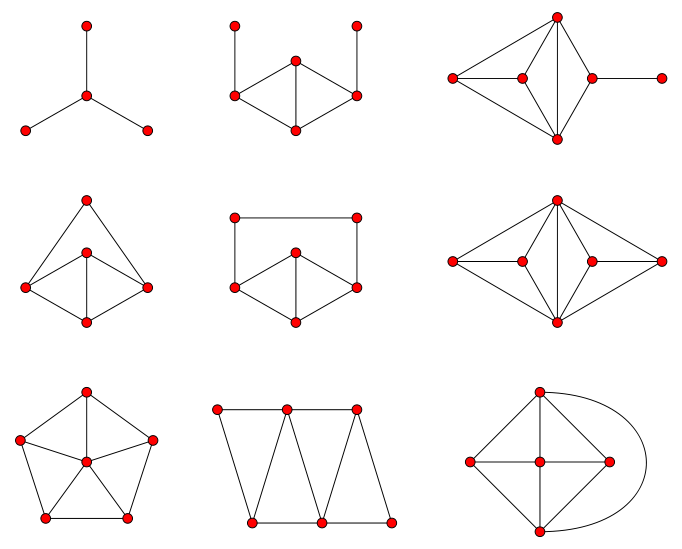}}
    \caption{Forbidden subgraphs for line graphs of graphs}
    \label{fig:l}   
\end{figure}

Likewise, the $FIS$ characterization of intersection graphs of  multigraphs were studied by Bermond and Meyer in 1973 in terms of the seven forbidden subgraphs [7]. Roussopoulos [43] and Lehot [26] gave a max (|E|, |V|) - time recognition algorithm for intersection graphs of graphs. In [43], the algorithm is based on the characterizations of intersection graphs due to Krausz whereas the algorithm in [26] is based on the characterizations of intersection graphs due to Rooij and Wilf. The algorithm from Degiorgi etl al. [9] is based on the local attributes. Naor et al. [39] proposed a parallel algorithm for line-to-root graph construction based on a divide-and-conquer scheme. Zverovich [52] gave a procedure for characterizing line graphs of a strict hereditary class of graphs in terms of FIS. This procedure is based on Beineke's characterization of line graphs of graphs [2] and Whitney's theorem on edge isomorphisms [51]. There are characterizations for line graphs of multigraphs with restricted multiplicities of edges (Tashkinov [45], bipartite graphs (Harary and Holzmann [15]) and bipartite multigraphs (Tyshkevich, Urbanovich and Zverovich [48]).

\section{Intersection graphs of k-uniform linear hypergraphs, $k = 3$.}

The situation changes if one takes $k = 3$ instead of $k = 2$. For $ k \geq 3$, the problem of characterizing $I(k)'s$ becomes complicated in terms of $F(I(k))'s$. A global characterization of representative graphs of k-uniform hypergraphs and that of k-uniform linear hypergraphs for an arbitrary k are given by Berge [4]. Lovasz [29] stated the problem of characterizing the class $I(3)$ intersection graphs of 3-uniform hypergraphs and noted that it cannot be characterized in terms of $FIS$. Bermond et al. [6], Gardner [11], Germa and Nickel [4] showed that this class cannot be characterized by a finite list of forbidden induced subgraphs $(FIS)$. It is also pointed out in [4] that any graph is a line graph for some linear hypergraph. For an arbitrary graph belonging to $I_{1} (k)$, the maximal cliques can be constructed in polynomial time (a polynomial in the vertex number, the power of the polynomial increases together with k) [27]. The recognition problem $G\in I(k)$ for fixed $k \geq 3 $ is NP-complete [27], [41].

The difficulty in finding a characterization for $I_{1}(3)$ is due to the fact that there are infinitely many forbidden induced subgraphs and their descriptions seem to be not simple. For an integer, $ m > 0$, consider a chain of $(m+2)$ diamond graphs such that the consecutive diamonds share vertices of degree two and add two pendant edges at every vertex of degree $2$ to get one of the families of minimal forbidden subgraphs of Naik et al. [37], [38] and [4] as shown below in figure 2 from [60]. Call this graph as $ {G_{1} (m)}$. Now let, ${\{G_{3}={G_{1} (m), m > 0}\}}$. It can be verified that $G_3  \in  F (I_1 (3))$. Additional examples can be found in [11].

\begin{figure}[!htb]
    \center{\includegraphics[width=0.5\textwidth] {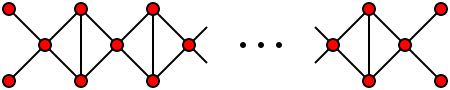}}
    \caption{A forbidden graph for intersection graphs of 3-uniform linear hypergraphs}
    \label{fig:2}   
\end{figure}

This does not rule out either the existence of polynomial time recognition $(ALG)$ or the possibility of a forbidden induced subgraph $(FIS)$ characterizations of k-uniform linear hypergraphs similar to Beineke's list of forbidden subgraphs. However, in [4] and [37], it is shown that $I_{1} (3)$ can be characterized by a finite list of forbidden induced subgraphs in the class of graphs whose vertex degrees are at least 69. In [37], and [38], authors also gave a global characterization of the members of the family $I_{1} (k)$, which is a generalization of a criterion of the intersection graphs due to Krausz. 

\begin{theorem}
$[37]$ If G is a graph, then $ G \in   I_{1} (k) $ if and only if in G there exists a set $C = \{C_1,\ldots,C_m \}$ of cliques of $sizes > 1$ such that the following condition holds. 

Every edge of G is in a unique click of C, and every vertex of G is in at most k cliques of C.
\end{theorem}

The following theorem $3.2$ similar to theorem $2.3$ of Rooij and Wilf is a generalization characterization of the family $I_{1} (k)$. 

\begin {theorem}
$[37]$ If G is a graph, then $G \in I_{1} (k)$ if and only if G has a set T of triangles satisfying the following two conditions:
\begin{enumerate}
    \item 
	If $abc$, and $abd$ are in T with $c \neq d$, then $ $cd$ \in E(G)$ and $acd$, $bcd$ are also in T;
	
	\item
	given any $(k+1)$ distinct edges of G, all having a vertex in common, at least two of these edges are in a triangle of G which is in T.
\end{enumerate}
\end{theorem}

The theorems $3.1$ and $3.2$ are used in proving the theorems $3.3$ and $4.2$ of section 4 for k > 3. The proofs of theorems $3.1$ and $3.2$ are similar to the proofs for $k = 2$ Harary [14]. 

A set of cliques of a graph G is a linear r-covering of G if G is the union of these cliques and any two cliques from the set have at most one common vertex and each vertex of G belongs to at most r cliques. A clique of size $ \geq    (k^2-k+2) $ is called k-large in $I_{1}(k)$. In $[23]$, $[27]$,  $[33-36]$, $[37-38]$, $[44]$ and $[53-57]$, it is noted that if G has a linear $k-$covering, then the covering must contain the set of all subgraphs of G induced by maximal $k-$ large cliques. 

\begin{theorem}
$[37]$ There is a finite family F(1) of forbidden graphs such that any graph G with minimum degree at least 69 belongs to $I_{1} (3)$ if and only if G has no induced subgraph isomorphic to a member of F(1).
\end{theorem}

Metelsky et al. [33] reduced the minimum degree bound from 69 to 19. Their proof is based on the theorems of Krausz, Beineke, and the Krausz characterization for the class $I_{1} (3)$ in $[37]$ and the properties of graph cliques from the class obtained by Levin et al. $[27]$. They constructed a finite family F(2) of forbidden graphs different from the family F(1) of the Theorem $3.3$. Their theorem states as follows.

\begin{theorem}
$[33]$ For a graph G with minimum degree $\geq 19$, the following two statements are equivalent.
\begin{enumerate}
    \item 
     $	G \in  I_{1} (3) $,
    \item 
    None of the graphs from the family F (2) is an induced subgraph of G.
\end{enumerate}
\end{theorem}

In $[33]$, the same authors gave the following polynomial algorithm $(ALG)$ for $I_{1} (3)$.

\begin{theorem}
$[33]$ There is a polynomial recognition algorithm to decide whether a graph with minimum degree $ \geq 19 $ belongs to $ I_{1} (3) $ .
\end{theorem}

To prove this, as noticed in $[37]$, the process is to construct a linear $3-$covering for the edges of G. 

Jacobson M. S. et al. [23] also gave a polynomial recognition algorithm $(ALG)$ for theorem $3.5$. The existence of the polynomial recognition algorithm follows from a simpler recursive characterization of graphs in $I_{1} (3)$ and relies on the fact that there is a polynomial time recognition algorithm for the members of $I_{1} (2)$.  

Denote by $ \delta_{alg} $ the minimal integer such that the problem $ G \in  I_{1} (3) $ is polynomially solvable in the class of graphs G with minimal degree $\geq  \delta_{alg} $ and by $ \delta_{fis} $ the minimal integer such that $ I_{1} (3) $ can be characterized by a finite list of forbidden induced subgraphs in the class of graphs G with minimal degree $\geq  \delta_{fis} $. The analysis of the constants $ \delta_{alg} $ and $ \delta_{fis} $ has splitted.

Matelsky et el. $[35]$ gave a polynomial algorithm solving the recognition problem for the graphs with bound on $\delta_{alg} \geq 13$. The complexity of the algorithm is $O(nm)$, where m and n are number of vertices and edges, respectively.The estimate on $\delta_{alg}$ is further improved: $6 \leq  \delta_{alg} \leq 11$ in  $[34]$.

Skums et al. $[44]$ further improved the bound on the minimum degree conditions. Their work deals mainly with the structural properties of the graphs from the class $I_{1} (3)$ connected with the geometry of the cliques. Their theorems are:

\begin{theorem}
$[44]$ There exists an algorithm with complexity $ O(nm) $ solving the recognition problem $ G \in  I_{1} (3) $ in the class of graphs G with minimum degree $\geq 10 $.
\end{theorem}

\begin{theorem}
$[44]$ If a graph G with minimum degree $ \geq 16 $ contains no graph from the set F(3) as an induced subgraph, then $G \in I_{1}(3) $.
\end{theorem}

The set F(3) of forbidden subgraphs $(FIS)$ is obviously different from the sets obtained previously on the minimum degree conditions.

Skums et al. $[44]$ proved $6 \leq \delta_{alg} $ $and$ $6 \leq \delta_{fis}$.

\section{Improvements on $Beineke's$ the nine forbidden subgraphs.}

Metelsky et al. $[33]$ reduced Beineke's nine forbidden subgraphs to six of them from figure 3. If the minimum degree of a graph G is $ \geq 5$, and G does not have the following induced subgraphs, then G is a line graph of a graph $[58]$. 

\begin{figure}[!htb]
    \center{\includegraphics[width=0.5\textwidth]
    {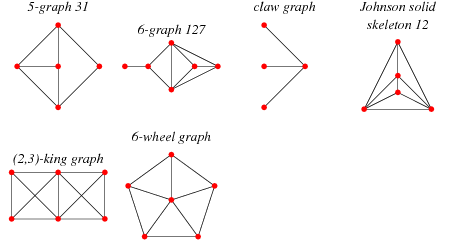}}
    \caption{Six forbidden subgraphs for line graphs with minimum degree five}
    \label{fig:3}   
\end{figure}

\section{Intersection graphs of k-uniform linear hypergraphs, $k > 3$.} 
Define inductively a family $G_{k}$, $k \geq 3$, as follows. $G_{3}$  is already defined above. 
For  k > 3, $G_k = \{\, G \mid  G \, $ by adding a pendent edge at every vertex of degree k in $G_1 $ where $G_1 \in G_k$ $ \} $. Clearly, $ G_{k} \in  F(I_{1} (k))$. Additional infinite families of forbidden graphs for $I_{l} (k)$  or $I(k)$ can be found in $[6]$, $[18]$, $[37-38]$. However, all these graphs have minimum degrees. Naik et al. raised a question in $[37]$ whether the theorem $3.3$ can be generalized for $k > 3$ or not. In the same paper of 1980, authors conjectured that the theorem $3.3$ cannot be generalized and later in 1997 this conjecture was proved by Matelsky et al. $[33]$. The proof of the conjecture is theorem $5.1$ stated below.

\begin{theorem}
$[33]$ For, $k > 3$, and an arbitrary constant $c$, the set of all graphs $G$ $\in I_{1} (k) $ with minimum degree $\geq c$  cannot be characterized by a finite list of forbidden induced subgraphs.
\end{theorem}

Metelsky et al. exhibited a new graph by taking the graph $G_{1}$ in $[37]$, $[38]$ as described above and pasting $(k-3)$ pairwise disjoint copies of the clique of size $s=max \{k^2-k+2,  c\}$  of $G_{1}$. Infinitely, many such graphs can be constructed. So, if $k > 3$, no such finite list characterizations exist for k-uniform linear hypergraphs, $k > 3$, no matter what lower bound is placed on the vertex degree.

It is proved in [19] that for any $k \geq 4$ and d > 0 the problem $G \in  I_{1} (k) $ remains NP-complete in the class of graphs G with the minimal vertex degree $\delta (G) \geq d $. 

To characterize $I_{1} (k)$, $k>3$, Naik, Rao, Shrikhande and Singhi developed a new method in $[38]$ based on the edge degree of a graph G. Edge degree of an edge e in G is the number of triangles in G containing the edge e and the minimum edge degree of G is the minimum over all the edges of G. Their characterization of  $I_{1} (k)$, $k\geq 3$ , based on the edge degree is stated in theorem 5.2.

\begin{theorem}
$[38]$, $[4]$ For a polynomial $f(k)=k^3-2k^2+1$, there exists a finite family $F(k)$ of forbidden graphs such that any graph $G$ with minimum edge degree $\geq  f(k)$ belongs to $I_{1} (k)$ if and only if $G$ has no induced subgraph isomorphic to a member of $F(k)$.
\end{theorem}

Jacobson et al.$[23]$ sharpened the bound on the polynomial $f(k)$ to $2k^2-3k+1$ by exhibiting a polynomial recognition algorithm $(ALG)$. Zverovich, I., in $[53]$ gave a forbidden graph characterization $(FIS)$ for $f(k) = 2k^2-3k+1$; their theorems are as follows.

\begin{theorem}
$[23]$ There is a polynomial algorithm $(ALG)$ to decide whether $ G \in I_{1}(k) $, $ k \geq 3 $ with edge degree  $ \geq   2k^2-3k+1 $.
\end{theorem}

\begin{theorem}
$[53]$ For $k \geq 3$, there is a finite set $F(k)$ of graphs such that a graph G with minimum edge degree $ \geq  2k^2-3k+1$ belongs to $I_{1} (k)$ if and only if G has no induced subgraph isomorphic to a member of  $F(k)$.
\end{theorem}

\begin{remark}
The complexity of recognizing intersection graphs of  $k-$uniform linear hypergraphs without any constraint on minimum vertex degree or minimum edge-degree is not known. \end{remark}

\section{Some open Problems and perspectives.}


\subsection{Naik et al. $[28]$ have the following problem connecting class of graphs $I_{1} (k)$ to the Eigen values of graphs. Eigen values of a graph are the Eigen values of its $(0, 1)$ adjacency matrix $[11]$. We will denote by $\alpha (G)$ the minimum Eigen value of graph G. For an arbitrary real number $\alpha$, we define, $E_{\alpha} $ \{\,$G \mid G$ is a graph with $\alpha(G) > \alpha$\,\}. It is evident that $E_{\alpha}$ is a hereditary family and that $I_{1} (k) \subseteq  E_{\alpha}$ for all $ \alpha < -k \leq -2$. They proposed a problem to describe the family F$(E_{\alpha})$ for all real number $\alpha$. Since $E_{\alpha}$ for $\alpha < -k $ is a larger family than $I_{1}(k)$, the set $F(E_{\alpha})$ may have simpler structures than that of $F(I_{1} (k))$. This is also suggested by Hoffman's theorem $[15]$,$[16]$,$[17]$. Hoffman's theorem, in fact, describes the families $E_{\alpha}$ in terms of families of $I_{l}(k)$.}

\subsection{Skums et al.[44] proposed to find the exact values for $\delta_{alg}$ and $\delta_{fis}$. They also proposed whether  $\delta_{alg} = \delta_{fis}$ or not. Furthermore, they made an interesting conjecture that $\delta_{alg}=9$.}

\subsection{Zverovich, I. $[53]$ proposed the following question on Theorem $5.4$. Is it possible to improve the bound $2k^2-3k+1$ on edge-degree? The guess is that the bound cannot be improved.}

\subsection{Zverovich, I. [55] proved that the class of graphs $C(k,l)$, $k \geq 0, l \geq 0$ branched out of Krausz partitions has $FIS$ characterizations. Author exhibited 14 $FISs$ for $C (3, 1)$, $C(k,l)$ is the class of all  graphs having a $l -$ bounded and $k -$ colorable.}

\subsection{In $[13]$ and $[19]$, the works on Krausz dimensions, computational complexities and other related topics branched out from the intersections graphs of $k- $ uniform linear hypergraphs are discussed.}

\subsection{Matelsky et. al. [36] studied the class $I_{2} (3)$ the intersection graphs of $k-$ uniform  hypergraphs, $k \le 3$ with multiplicity at most 2 (non-linear) and characterized them by a means of finite list of forbidden induced subgraphs (FIS) in the class of threshold graphs $[8]$,$[10]$. $ O(n) - $ time algorithm is given.}








\end{document}